\documentclass[12pt,oneside,english]{amsart}
\usepackage[latin1]{inputenc}
\usepackage{amssymb}

\makeatletter
\newtheorem{theorem}{Theorem}
\newtheorem{lemma}{Lemma}
\newtheorem{example}{Example}
\newtheorem{definition}{Definition}
\newtheorem{criterion}{Criterion}

\newtheorem{proposition}{Proposition}

\newtheorem{conjecture}{Conjecture}

\newtheorem{corollary}{Corollary}

\numberwithin{equation}{section}
\usepackage{babel}

\makeatother
\begin{document}
\baselineskip=17pt

\title[Three theorems on twin primes]{Three theorems on twin primes }

\author{Vladimir Shevelev}
\address{Departments of Mathematics \\Ben-Gurion University of the
 Negev\\Beer-Sheva 84105, Israel. e-mail:shevelev@bgu.ac.il}


\begin{abstract}
For earlier considered sequence: $c(1)=2$ and for $n\geq 2,$
$$c(n)=c(n-1)+\begin{cases}\gcd (n, \enskip c(n-1)),\enskip if \;\; n\enskip is\enskip even
\\\gcd (n-2,\enskip c(n-1)),\enskip if\;\; n\enskip is\enskip odd\end{cases} $$
we prove theorems of its connection with twin primes. We also give a sufficient condition for the infinity of twin primes and pose several new conjectures; among them we propose a very simple conjectural algorithm of constructing a pair $(p,\enskip p+2)$ of twin primes over arbitrary given integer $m\geq4$ such that $p+2\geq m.$
\end{abstract}
\maketitle

\section{Introduction }

In [2] we posed the following conjecture
\begin{conjecture}\label{1}
Let $c(1)=2$ and for $n\geq 2,$
$$c(n)=c(n-1)+\begin{cases}\gcd (n, \enskip c(n-1)),\enskip if \;\; n\enskip is\enskip even
\\\gcd (n-2,\enskip c(n-1)),\enskip if\;\; n\enskip is\enskip odd\end{cases}. $$
 Then every record (more than 3) of the values of difference $c(n)-c(n-1)$ is greater of twin primes.
 \end{conjecture}
 The first such records are (cf. sequence A166945 in [4])
 \begin{equation}\label{1.1}
  7,13,43,139,313,661,1321,2659,5459,10891,22039,...
  \end{equation}
  Our observations of the behavior of sequence $\{c(n)\}$ are the following:\newline\newline

  1) In some sequence of arguments $\{m_i\}$ we have $\frac {c(m_i)} {m_i}=2.$ These values of arguments we call \slshape the fundamental points. \upshape The least fundamental point is $m_1=2.$\newline
  2)For every two adjacent fundamental points $m_j<m_{j+1},$ we have $m_{j+1}\geq 2m_j.$\newline
  3) For $i\geq2,$ the numbers $m_i\mp1$ are twin primes (and, consequently, $m_i\equiv0\pmod6$).\newline
  4) In points $m_i+3$ we have $c(m_i+3)-c(m_i+2)=m_i+1.$ These increments we call \slshape the main increments\upshape \enskip of sequence $\{c(n)\},$ while other nontrivial (i.e.more than 1) increments we call \slshape the minor increments.\upshape\newline
  5)For $i\geq2,$ denote $h_i$ the number of minor increments between adjacent fundamental points $m_i$ and $m_{i+1}$ and $T_i$ the sum of these increments. Then $T_i\equiv h_i\pmod6.$\newline
  6)For $i\geq3,$ the minor increments between adjacent fundamental points $m_i$ and $m_{i+1}$ could occur only before $m_{i+1}-\sqrt{m_{i+1}-1}-4.$ \upshape\newline\newline
  \indent The aim of this paper is to show that the validity of all these observations follow only from 6).
  \begin{theorem}\label{1}
If observation 6) is true then observation 1)-5) are true as well.
 \end{theorem}
  \begin{corollary}\label{1}
 If $1)$ observation 6) is true and $2)$ the sequence $\{{c}(n)\}$ contains infinitely many fundamental points, then there exist infinitely many twin primes.
\end{corollary}
 Besides, in connection with Conjecture 1 we think that
 \begin{conjecture}\label{2}
  For $n\geq1,$ the main and only main increments are the record differences $c(n)-c(n-1).$
 \end{conjecture}

\section{Proof of Theorem 1}
Note that $c(12)=24$ and the numbers  $11,13$ are primes. We use induction. Suppose $n_1\geq15$ is a number of the form 6l+3 (for $n_1<15$ the all observations are verified directly).
Let $n_1-3$ is a fundamental point: $$c(n_1-3)=2n_1-6$$ and for $n:=n_1-3,\enskip n\mp1$ are twin primes. Since $n_1-2$ is odd, then we have
$$c(n_1-2)=2n_1-5.$$
Further, since $\gcd(n_1-1,\enskip2n_1-5)=\gcd(n_1-1,\enskip3)=1,$ then
$$c(n_1-1)=2n_1-4,$$
Since  $$\gcd (n_1-2,\enskip 2n_1-4)=n_1-2$$ then we have a main increment such that
 \begin{equation}\label{2.1}
c(n_1)=3n_1-6.
 \end{equation}

Here we distinguish two cases:\newline
\bfseries A \mdseries) Up to the following fundamental point there are only trivial increments. The inductive step in this case we formulate as the following.
 \begin{theorem}\label{2}
If $6\leq m_j<m_{j+1}$ are adjacent fundamental points with only nontrivial increment between them which is a main increment, then \newpage
i)\enskip$m_{j+1}=2m_j;$\newline
ii)\enskip If $m_j\mp1$ are twin primes, then $m_{j+1}\mp1$ are twin primes as well.
 \end{theorem}
  Thus, if observation 2) is true, then to every pair of the adjacent fundamental points with the only main increment between them corresponds a quadruple of primes of the form
 $$p,\enskip p+2,\enskip 2p+1,\enskip 2p+3.$$

 \begin{example} Consider the adjacent fundamental points $n_1=660$ and $n_2=1320.$ Since $n_2=2n_1,$ then between them there is no any miner increment. We have
 $$c(660)=1320,\enskip c(661)=1321,\enskip c(662)=1322,\enskip c(663)=1983$$ and
 $$c(663)-c(662)=661.$$
 $$c(1320)=2640,\enskip c(1321)=2641,\enskip c(1322)=2642,\enskip c(1323)=3963$$ and
 $$c(1323)-c(1322)=1321.$$
 Here we have two pairs of twins:
 $$p=659, \enskip p+2=661, \enskip 2p+1=1319, \enskip 2p+3=1321.$$
\end{example}
\bfseries Inductive step in case A \mdseries)\newline Continuing (2.1), we have
$$c(n_1+1)=3n_1-5,$$
$$c(n_1+2)=3n_1-4,$$
$$...$$
$$c(2n_1-6)=4n_1-12,$$
(It is the second fundamental point in the inductive step)
$$c(2n_1-5)=4n_1-11,$$
$$c(2n_1-4)=4n_1-10,$$
Since
 $$\gcd (2n_1-5,\enskip 4n_1-10)=2n_1-5$$ then, denoting $n_2=2n_1-3,$ we have
\begin{equation}\label{2.2}
c(2n_2)=3n_2-6.
 \end{equation}
 Note that, since $n_1=6l+3,$ then $n_2=6l_1+3,$ where $l_1=2l.$\newpage

\indent Furthermore, from the run of formulas (2.2) we find for $5\leq j\leq\frac{n_1-2} {2}$
$$c(2n_1-2j-1)=4n_1-2j-7,$$
$$c(2n_1-2j)=4n_1-2j-6.$$
This means that
$$\gcd(2n_1-2j,\enskip 4n_1-2j-7)=1,\enskip i.e.\enskip \gcd(2j-7,\enskip 2n_1-7)=1.  $$
Note that, for the considered values of $n_1$ we have $2\frac{n_1-2} {2}-7\geq\sqrt{2n_1-7},$ then $2n_1-7=n_2-4$ is prime.\newline
\indent On the other hand,
$$c(2n_1-2j)=4n_1-2j-6,$$
$$c(2n_1-2j+1)=4n_1-2j-5.$$
Thus, for $5\leq j\leq\frac{n_1-1} {2},$
$$\gcd(2n_1-2j-1,\enskip 4n_1-2j-6)=1,\enskip i.e.\enskip \gcd(j-2,\enskip 2n_1-5)=1.$$

Here, for the considered values of $n_1$ we also have $\frac{n_1-5} {2}\geq\sqrt{2n_1-5},$ then $2n_1-5=n_2-2$ is prime as well. This completes the inductive step in case \bfseries A \mdseries).
 If, in addition, to note that $n_1-3$ and $n_2-3$ are the two adjacent fundamental points, then we get a proof of Theorem 2.$\blacksquare$

 \bfseries B \mdseries) Up to the following fundamental point we have some minor increments.

   The inductive step we formulate as following.
   \begin{theorem}\label{3}
Let observation 6) be true. If $6\leq m_i<m_{i+1}$ are adjacent fundamental points with a finite number of minor increments between them, then \newline
i)\enskip$m_{i+1}\geq2m_i;$\newline
ii)\enskip If $m_i\mp1$ are twin primes, then $m_{i+1}\mp1$ are twin primes as well.
 \end{theorem}

 Thus the observation 2) will be proved in frameworks of the induction.\newline

 \bfseries Inductive step in case B \mdseries)

 Let in the points $n_1+l_j\enskip j=1,...,h,$ before the second fundamental point we have the minor increments $t_j,\enskip j=1,...,h.$ We have ( starting with the first fundamental point $n_1-3$)
$$c(n_1-3)=2n_1-6$$
$$c(n_1-2)=2n_1-5,$$
$$c(n_1-1)=2n_1-4,$$
$$c(n_1)=3n_1-6,$$\newpage
$$c(n_1+1)=3n_1-5,$$
$$...$$

$$c(n_1+l_1-1)=3n_1+l_1-7.$$

\begin{equation}\label{2.3}
c(n_1+l_1)=3n_1+l_1+t_1-7,
\end{equation}

$$c(n_1+l_1+1)=3n_1+l_1+t_1-6,$$

$$...$$

$$c(n_1+l_2-1)=3n_1+l_2+t_1-8,$$

\begin{equation}\label{2.4}
c(n_1+l_2)=3n_1+l_2+t_1+t_2-8,
\end{equation}
$$...$$
$$c(n_1+l_h-1)=3n_1+l_h+t_1+...+t_{h-1}-h-7,$$
\begin{equation}\label{2.5}
c(n_1+l_h)=3n_1+l_h+t_1+...+t_h-h-6,
\end{equation}
$$c(n_1+l_h+1)=3n_1+l_h+t_1+...+t_h-h-5,$$
$$...$$
\begin{equation}\label{2.6}
c(2n_1+T_h-h-6)=4n_1+2T_h-2h-12,
\end{equation}
where
\begin{equation}\label{2.7}
T_h=t_1+...+t_h.
\end{equation}
(thus $2n_1+T_h-h-6$ is the second fundamental point in the inductive step)
$$c(2n_1+T_h-h-5)=4n_1+2T_h-2h-11,$$
$$c(2n_1+T_h-h-4)=4n_1+2T_h-2h-10.$$ Here we need a lemma.
\begin{lemma}
$T_h-h $ is even.
\end{lemma}
\bfseries Proof.  \mdseries We use the induction over $h\geq1.$ If $l_1$ is even, then, by (2.3),
$$\gcd(n_1+l_1-2,\enskip 3n_1+l_1-7)=t_1$$
and $t_1$ divides $2n_1-5.$ Analogously, if $l_1$ is odd, then $t_1$ divides $2n_1-7.$ \newpage Thus $T_1-1=t_1-1$ is even. Suppose that $T_{h-1}-(h-1)$ is even. Then, by (2.5) in the case of $l_h$ is even, we have
$$\gcd(n_1+l_h-2,\enskip 3n_1+l_h+T_{h-1}-(h-1)-7)=t_h$$
and $t_h$ divides, by the inductive supposition, an odd number $2n_1+T_{h-1}-(h-1)-5.$ Analogously, if $l_h$ is odd, then $t_h$ divides odd number $2n_1+T_{h-1}-(h-1)-7.$ Thus $t_h$ is odd and we conclude that $T_h-h=T_{h-1}-(h-1)+t_h-1$ is even.$\blacksquare$\newline
 Therefore, we have
 $$\gcd (2n_1+T_h-h-5,\enskip4n_1+2T_h-2h-10)=2n_1+T_h-h-5,$$
 and in the point $n_2:=2n_1+T_h-h-3$ we have the second main increment (in framework of the inductive step). Thus
\begin{equation}\label{2.8}
c(2n_1+T_h-h-3)=6n_1+3T_h-3h-15.
\end{equation}
 Note that, for $n\geq2,$ we have $c(n)\equiv n\pmod2.$ Therefore, $T_h\geq3h$ and for the second fundamental point $n_2-3=2n_1+T_h-h-6$ we find

$$n_2-3\geq2(n_1-3)+2h.$$

 By the induction (with Theorem 2), \slshape this proves observation 2).\upshape\newline
 \indent Now, in order to finish the induction, we prove the primality of numbers $n_2-4=2n_1+T_h-h-7$ and  $n_2-2=2n_1+T_h-h-5.$\newline
\indent From the run of formulas (2.5)-(2.6) for $5\leq j\leq\frac{n_1+T_h-h-l_h-3} {2}$
(unfortunately,we cannot cross the upper boundary of the last miner increment) we find
$$c((2n_1+T_h-h-4)-(2j-1))=4n_1+2T_h-2h-2j-9,$$
$$c(2n_1+T_h-h-2j-2)=4n_1+2T_h-2h-2j-8.$$
This means that
$$\gcd(2n_1+T_h-h-2j-2,\enskip 4n_1+2T_h-2h-2j-9)=1,$$ i.e.\enskip $$\gcd(2j-7,\enskip 2n_1+T_h-h-7)=1.  $$
For the most possible $j=\frac{n_1+T_h-h-l_h-3} {2}$ we should have
$$2j-7=n_1+T_h-h-l_h-10\geq\sqrt{2n_1+T_h-h-7},$$
or, since $n_2=2n_1+T_h-h-3,$ then we should have $n_2-n_1-l_h-7\geq\sqrt{2n_1+T_h-h-7},$
 i.e.
$$n_1+l_h\leq n_2-\sqrt{n_2-4}-7.$$

 This condition is equivalent to the observation 6) which is written in\newpage terms of the fundamental points $m_i=n_i-3.$
Thus from observation 6) we indeed obtain the primality of $n_2-4=2n_1+T_h-h-7.$ \newline
\indent What is left-to prove the primality of $n_2-2=2n_1+T_h-h-5.$ We do it in the next section \slshape without supposition of the validity of observation 6).\upshape

\section{Completion of proof of Theorem 1: proof of the primality of $2n_1+T_h-h-5$ independently on observation 6)}
It is interesting that, using the Rowland's method [1], we are able to get  the primality of $2n_1+T_h-h-5$ without unproved observation 6). This gives an additional hope to convert the observations 1)-6) into the absolute statements.

 Denote
 \begin{equation}\label{3.1}
  n_1^*:=n_1+l_h.
  \end{equation}
By (2.5),
\begin{equation}\label{3.2}
  c(n_1^*)= 3n_1+l_h+T_h-h-6=3n_1^*-2l_h+T_h-h-6=3n_1^*+u,
 \end{equation}
 where
 \begin{equation}\label{3.3}
  u=T_h-h-2l_h-6\equiv0\pmod2,
  \end{equation}
 and
  \begin{equation}\label{3.4}
  c(n_1^*+i-1)= 3n_1^*+i+u-1, \enskip i\leq k,
 \end{equation}
 where $k$ is the smallest positive integer such that the point $3n_1^*+k+u-1$ is point of a nontrivial increment.
Put $h(n)=c(n)-c(n-1),$ such that
  \begin{equation}\label{3.5}
 h(n)=\begin{cases}\gcd (n, \enskip c(n-1)), \enskip if \;\; n\enskip is\enskip even,
\\\gcd (n-2,\enskip c(n-1)), \enskip if\;\; n\enskip is\enskip odd,\end{cases}
 \end{equation}
 then
 $$h(n_1^*+i)=\begin{cases}\gcd (n_1^*+i,\enskip c(n_1^*+i-1)), \enskip if \;\; n_1^*+i\enskip is\enskip even,
\\\gcd (n_1^*+i-2,\enskip c(n_1^*+i-1)), \enskip if\;\; n^*_1+i\enskip is\enskip odd.\end{cases}  $$
Put
\begin{equation}\label{3.6}
 \delta=\delta(m)=\begin{cases}0, \enskip if \;\; m\enskip is\enskip even,
\\2, \enskip if\;\; m\enskip is\enskip odd.\end{cases}
 \end{equation}
Thus, $ h(n_1^*+i)$ divides both $n_1^*+i-\delta(n_1^*+i)$ and $3n_1^*+i+u-1$ and also divides both their difference
\begin{equation}\label{3.7}
 2n_1^*+u-1+\delta(n_1+i)
 \end{equation}
and
\begin{equation}\label{3.8}
 3(n_1^*+i-\delta(n_1^*+i))-(3n_1^*+i+u-1)=2i-u+1-3\delta(n_1^*+i).
 \end{equation}\newpage
 Let $q$ is the smallest prime divisor of
 \begin{equation}\label{3.9}
 c(n_1^*)-n_1^*+\delta(n_2)-1=(by \enskip (3.2))\enskip 2n_1^*+u+\delta(n_2)-1.
 \end{equation}
Note that, since $u$ is even, then $q$ is odd. Let us prove that
\begin{equation}\label{3.10}
 k\geq\frac{1} {2}(q+u-1+3\delta(n_1^*+k)).
 \end{equation} In view of (3.7), for $ i\leq k,$ the number $ h(n_1+i)$  divides $2n_1^*+u+\delta(n_1^*+i)-1.$ Therefore, for $i=k,$ we see that
\begin{equation}\label{3.11}
h(n_1^*+k)\geq q.
 \end{equation}
 Since , by (3.4), $h(n_1^*+k)$ divides $2k-u+1-3\delta(n_1^*+k),$ then, using (3.11), we find
$$q\leq h(n_1^*+k)\leq 2k-u+1-3\delta(n_1^*+k)$$
and (3.10) follows.\newline
Now show that also
 \begin{equation}\label{3.12}
 k\leq\frac{1} {2}(q+u-1+3\delta(n_1^*+k)).
 \end{equation}
 By the definition of $k,$ for $1\leq i<k,$ we have $h(n_1^*+i)=1,$ and, using (3.10), we conclude that at least for $1\leq i<\frac{1} {2}(q+u-1+3\delta(n_1^*+i))$ we have $h(n_1^*+i)=1.$ Show that $i=\frac {1} {2}(q+u-1+3\delta(n_1^*+k))$ produces a nontrivial
 $\gcd.$ Indeed, according to (3.5), we have
$$h(n_1^*+ \frac{1} {2}(q+u-1+3\delta(n_1^*+k)))=$$
$$\gcd(n_1^*-\delta(n_1^*+k)+\frac{1} {2}(q+u-1+3\delta(n_1^*+k)),$$
$$ 3n_1^*+u-1+\frac{1} {2}(q+u-1+3\delta(n_1^*+k)))=$$

 $$\gcd(\frac{1} {2}((2n_1^*+u+\delta(n_1^*+k)-1)+q),$$
 \begin{equation}\label{3.13}
  \frac{1} {2}(3(2n_1^*+u+\delta(n_1^*+k)-1)+q)).
 \end{equation}
 From (3.10) and (3.13) it follows that $q$ divides both of arguments of $\gcd.$ Therefore,
 $$h(n_1^*+ \frac{1} {2}(q+u-1+3\delta(n_1^*+k)))\geq q\geq3. $$
 Thus, by the definition of number $n_2,$ we have
 $$ k=n_2-n_1^*\leq\frac{1} {2}(q+u-1+3\delta(n_1^*+k)).$$
 Therefore, \begin{equation}\label{3.14}
 k=\frac{1} {2}(q+u-1+3\delta(n_1^*+k)).
 \end{equation}\newpage
 On the other hand, according to (3.8), $h(n_1^*+k),$ divides $2k-u+1-3\delta(n_1^*+k),$ or, taking into account (3.14), divides $q.$ Therefore,
 \begin{equation}\label{3.15}
  h(n_2)=h(n_1^*+k)=q.
 \end{equation}
 According to (3.14)-(3.15), we have
  \begin{equation}\label{3.16}
 h(n_2)=q=2k-u-3\delta(n_2)+1.
 \end{equation}
 Nevertheless, by (2.8), $n_2=2n_1+T_h-h-3\equiv1\pmod2$ and, by (3.3) $u=T_h-h-2l_h-6.$ Therefore,
  \begin{equation}\label{3.17}
 q=2k-u-5=2(n_2-n_1^*)-u-5=2n_1+T_h-h-5.
 \end{equation}
 Thus $2n_1+T_h-h-5$ is prime. This completes proof of Theorem 1.$\blacksquare$\newline\newline
\begin{corollary}\label{2}
 If $p_1<p_2$ are consecutive greater of twin primes giving by Theorem 1, then $p_2\geq 2p_1-1.$
\end{corollary}
\bfseries Proof.  \mdseries Since (see already proved observation 2) $n_2-3\geq2(n_1-3),$ then $q=n_2-2\geq2(n_1-2)-1,$ where, by the inductive supposition, $ n_1-2$ is greater of twin primes. Now the corollary follows in the frameworks of the induction.$\blacksquare$
\begin{corollary}\label{3}
 $$T_h\equiv h\pmod6.$$
 \end{corollary}
\bfseries Proof.  \mdseries  The corollary immediately follows from the well known fact that the half-sum of twin primes not less than 5 is a multiple of 6. Therefore,  $2n_1+T_h-h-6\equiv0\pmod6.$ Since, by the condition, $2n_1\equiv0\pmod6,$ then we obtain the corollary.$\blacksquare$\newline
\indent Now the observation 5) follows in the frameworks of the induction. The same we can say about observation 4).

\indent The observed weak excesses of the exact estimate of Corollary 2 indicate to the smallness of $T_h$ and confirm, by Theorem 1, Conjecture 1.
\section{Estimates of ratios $c(n)/n$ and stronger conjecture}
From the construction of Section 2 it easily follows that only in the fundamental points of the considered sequence we have $\gamma(n):=c(n)/n=2.$ Moreover, only in two points following after every fundamental point we have the values of $\gamma(n)$ less than 2. Namely, if $n$ is a fundamental point, then in the point $\nu=n+1$ we have $\gamma(\nu)=2-\frac{1} {\nu}$ and in the point $\mu=n+2$ we have $\gamma(\mu)=2-\frac{2} {\mu}.$ On the other hand,
using induction, it is easy to\newpage prove that
$$\frac {c(n)} {n}\leq \begin{cases}3,\enskip if \;\; n\enskip is\enskip even,
\\3-\frac{6} {n},\enskip if\;\; n\enskip is\enskip odd.\end{cases} $$
Indeed, let
$$c(n-1)\leq \begin{cases}3(n-1),\enskip if \;\; n\enskip is\enskip even,
\\3(n-1)-6,\enskip if\;\; n\enskip is\enskip odd.\end{cases} $$

Since $$h(n)=c(n)-c(n-1)|\begin{cases}n\enskip and \enskip c(n-1),\enskip if \;\; n\enskip is\enskip even,
\\n-2\enskip and \enskip c(n-1),\enskip if\;\; n\enskip is\enskip odd,\end{cases} $$
then
 $$h(n)\leq\begin{cases}3n-c(n-1),\enskip if \;\; n\enskip is\enskip even,
\\3n-6-c(n-1),\enskip if\;\; n\enskip is\enskip odd,\end{cases} $$
and
$$c(n)=c(n-1)+h(n)\leq \begin{cases}3(n-1),\enskip if \;\; n\enskip is\enskip even,
\\3(n-1)-6,\enskip if\;\; n\enskip is\enskip odd.\end{cases} $$
Thus we proved the following estimates.

\begin{proposition}
\begin{equation}\label{4.1}
2-\frac {2} {n-1}\leq\gamma(n)\leq \begin{cases}3,\enskip if \;\; n\enskip is\enskip even,
\\3-\frac{6} {n},\enskip if\;\; n\enskip is\enskip odd.\end{cases}
 \end{equation}
\end{proposition}
In points $n$ of the main increments we have $\gamma(n)=3-\frac{6} {n}.$ The first terms of the sequence $\{\beta_j\}$ for which  $\gamma(\beta_j)=3$ are:
$$18,\enskip 20,\enskip 66,\enskip 150,\enskip 156,\enskip 1326,\enskip 10904,\enskip 10908,\enskip 10910,... $$\newline

It is easy to see that observation 6) one can replace by, e.g., the observation that, for every $i\geq3,$ \slshape in the maximal point $\rho^{(i)}$ of a nontrivial increment before fundamental point $m_i$ we have \upshape
\begin{equation}\label{4.2}
\gamma(\rho^{(i)})\geq 2.5\enskip.
\end{equation}
Indeed, putting in (2.5) $n_1:=n_{i-1}=m_{i-1}+3,\enskip n_2:=n_{i}=m_{i}+3$ and $n_1+l_h:=\rho^{(i)},$ such that, by (2.7)(see the second fundamental point of the inductive process), $T_h-h-6:=m_{i}-2n_{i-1}$ we, by the supposition, have
\begin{equation}\label{4.3}
\gamma(\rho^{(i)})=\frac {\rho^{(i)}+m_{i}} {\rho^{(i)}}\geq 2.5\enskip.
\end{equation}
Thus
\begin{equation}\label{4.4}
 \rho^{(i)}\leq \frac{2} {3}m_{i}.
\end{equation}\newpage
This means that the distance between $\rho^{(i)}$ and $m_{i}$ is not less than $m_{i}/3.$ Since we have $x/3>\sqrt{x-1}+4,$ for $x\geq30,$ then observation 6) follows for $m_{i+1}\geq30.$  \newline
\indent  Our stronger conjecture is the following.
\begin{conjecture}\label{3}
Let $m_{i-1}<m_{i}$ be adjacent fundamental points. Let $\rho^{(i)}$ be the maximal point of a nontrivial increment before $m_{i}.$ Then
 \begin{equation}\label{4.5}
\lim_{i\rightarrow\infty}\frac {\rho^{(i)}} {m_i}=\frac {1} {2}.
\end{equation}
 \end{conjecture}

\section{A sufficient condition for the infinity of twin primes}

Put
\begin{equation}\label{5.1}
 \rho^{(i)}=\lambda_i n_{i-1}, \enskip\lambda_i\geq1.
 \end{equation}

\begin{conjecture}\label{4}
For every $i\geq6,\enskip \lambda_i\leq5/4.$
 \end{conjecture}
\begin{theorem}
If Conjecture 4 is true, then we have infinitely many twin primes.
\end{theorem}
\bfseries Proof.  \mdseries Since
\begin{equation}\label{5.2}
 c(\lambda_i n_{i-1})=\gamma(\rho^{(i)})\lambda_i n_{i-1},
 \end{equation}
 then the distance $r_i$ between $\rho^{(i)}$ and $m_{i}$ is defined by the equation
 \begin{equation}\label{5.3}
\frac{\gamma(\rho^{(i)})\lambda_i n_{i-1}+r_i} {\lambda_i n_{i-1}+r_i}=2.
 \end{equation}
Thus
\begin{equation}\label{5.4}
r_i=\gamma(\rho^{(i)})\lambda_i n_{i-1}-2\lambda_i n_{i-1}
 \end{equation}
and we have
 \begin{equation}\label{5.5}
m_{i}=\lambda_i n_{i-1}+r_i=\lambda_i n_{i-1}(\gamma(\rho^{(i)})-1).
 \end{equation}
Put
\begin{equation}\label{5.6}
m_{i}=\lambda_i n_{i-1}(\gamma(\rho^{(i)})-1)=(2+\mu_i)n_{i-1}.
 \end{equation}
Since, by Theorem 2, which was proved independently from observation 6), we have
$$2\leq\frac {m_{i}} {m_{i-1}}=\frac {m_{i}} {n_{i-1}-3},$$
then $m_{i}\geq2n_{i-1}-6$ and, by (5.6),
\begin{equation}\label{5.7}
2+\mu_i=\lambda_i (\gamma(\rho^{(i)})-1)\geq2-\frac{6} {n_{i-1}}.
 \end{equation}
Furthermore, by the condition, $\lambda_i\leq5/4.$  Therefore, we have\newpage
$$5/4\geq \lambda_i\geq(2-\frac {6}{n_{i-1}})\frac{1}{\gamma(\rho^{(i)})-1}. $$
Note that, for $i\geq6,$ we have $n_{i-1}\geq141.$
Therefore,

\begin{equation}\label{5.8}
\gamma(\rho^{(i)})\geq1+\frac{4} {5}(2-\frac {6}{n_{i-1}})\geq5/2.
 \end{equation}
 By (4.3)-(4.4), this means that observation 6) follows and the numbers $m_i\mp1$ are twin primes.\newline
 On the other hand, by (5.5) and Proposition 1, we have
 \begin{equation}\label{5.9}
 m_{i}=\lambda n_{i-1}(\gamma(\rho^{(i)})-1)\leq2\lambda_i n_{i-1}\leq2.5 n_{i-1}.
\end{equation}
 The latter inequality ensures the infinity of the fundamental points of the considered sequence and, consequently, the infinity of twin primes.$\blacksquare$\newline\newline
 \indent Moreover, if Conjecture 4 is true, then verifying a finite set of integers beginning with $n=2,$ \enskip from Theorem 4 we obtain that:\newline\newline\indent\slshape Between $n\geq2$ and $3n$ we have at least one pair of twin primes.\upshape\newline\newline
\indent Note that, the first real values of $\lambda_i=\frac{\rho^{(i)})} {n_{i-1}},\enskip i\geq6$ are:

$$\frac{156} {141}=1.106...;\enskip\frac{348} {315}=1.104...;\enskip \frac{661} {661}=1.000...;$$
 $$\enskip \frac{1339} {1323}=1.012...;\enskip\frac{2712} {2661}=1.019...;\enskip \frac{5496} {5421}=1.013...\enskip,...$$\newline\newline
 \indent Note that if the last miner increment $\rho^{(i)}$ after the point of the main increment $n_{i-1}$ is known, then the following miner increment is
 \begin{equation}\label{5.10}
  n_i=c(\rho^{(i)})-\rho^{(i)}.
 \end{equation}
 It easy follows from (2.5)-(2.6).

\section{To every integer $m\geq4$ corresponds a pair of twin primes $(p,\enskip p+2)$ such that $p+2\geq m$}
Given $m\geq4,$ we give a very simple rule to calculate a pair of twin primes (p,\enskip p+2) such that $p+2\geq m.$
Although till now we are able to prove a private case of this rule, we absolutely do not doubt that it is always true!
For every positive integer $m,$ consider the following sequence:\newpage
$$c^{(m)}(1)=m;\enskip for\enskip n\geq2,$$
 \begin{equation}\label{6.1}
c^{(m)}(n)=c^{(m)}(n-1)+\begin{cases}\gcd (n, \enskip c^{(m)}(n-1)),\enskip if \;\; n\enskip is\enskip even
\\\gcd (n-2,\enskip c^{(m)}(n-1)),\enskip if\;\; n\enskip is\enskip odd.\end{cases}
\end{equation}
Thus for every $m$ this sequence has the the same formula that the considered one but another initial condition.
Our astonishing observation is the following.
\begin{conjecture}\label{5}
Let $n^*,$ where $n^*=n^*(m),$ be point of the last nontrivial increment of $\{c^{(m)}(n)\}$ on the set $N_m=\{1,...,m\}$ and $n^*=1,$ if there is not any nontrivial increment on $N_m.$
Then numbers $c^{(m)}(n^*)-n^*\mp1$ are twin primes.
 \end{conjecture}
Evidently, $c^{(m)}(n^*)-n^*+1\geq m$ and the equality holds if and only if $n^*=1.$
 \begin{example}\label{2}
Let $m=20.$ Then $n^*=12$ and $c^{(m)}(n^*)=42.$ Thus numbers $42-12\mp1$ are twin primes $(29,\enskip31).$
 \end{example}
  \begin{example}\label{3}
Let $m=577.$ Then $n^*=156$ and $c^{(m)}(n^*)=1038.$ Thus numbers $1038-156\mp1$ are twin primes $(881,\enskip883).$
 \end{example}
 \begin{example}\label{4}
Let $m=3000.$ Then $n^*=2$ and $c^{(m)}(n^*)=3002.$ Thus numbers $3002-2\mp1$ are twin primes $(2999,\enskip3001).$
 \end{example}
 The case of $n^*=1$ we formulate as the following criterion.
 \begin{criterion}\label{1} A positive integer $m>3$ is a greater of twin primes if and only if the points $1,...,m$ are points of trivial increments of sequence $\{c^{(m)}(n)\}.$
\end{criterion}
\bfseries Proof.  \mdseries
By the condition,
$$c^{(m)}(1+i)=m+i,\enskip c^{(m)}(2+i)=m+i+1,$$
Therefore, if $i$ is even, then
$$\gcd(2+i,\enskip m+i)=1,$$
or
$$\gcd(m-2,\enskip i+2)=1.$$
If $i$ is odd, then
$$\gcd(i,\enskip m+i)=1,$$
or
$$\gcd(m,\enskip i)=1.$$\newpage
Since $i$ is arbitrary from $N_m,$ then both of numbers $m-2,\enskip m$ are primes.
The converse statement is also evident.$\blacksquare$\newline
\section{A theorem on twin primes which is independent on observation of type 6)}
Here we present a new sequence $\{a(n)\}$ with the quite analogous definition of fundamental and miner points for which Corollary 1 is true in a stronger formulation. Using a construction close to those ones that we considered in [3], consider the sequence defined as the following:

 $a(180)=360$ and for $n\geq 181,$
\begin{equation}\label{7.1}
a(n)=\begin{cases}a(n-1)+1,\enskip if \;\;\gcd(n+(-1)^n-1,\enskip a(n-1))=1;
\\3(n-2)\;\; otherwise\end{cases}.
\end{equation}

\begin{definition} A point $m_i$ is called \upshape a fundamental point of sequence (7.1),\slshape \enskip if it has the form $m_i=6t$ and $a(m_i)=2m_i.$ The increments in the points $m_i+3$ we call the \upshape main increments.\slshape \enskip Other nontrivial increments we call \upshape miner increments.
\end{definition}
The first fundamental point of sequence (7.1) is $m_1=180.$

\begin{theorem}\label{5}
If the sequence $\{{a}(n)\}$ contains infinitely many fundamental points, then there exist infinitely many twin primes.
 \end{theorem}
 \bfseries Proof.  \mdseries We use induction. Note that numbers $m_1\mp1$ are twin primes: 179 and 181. Suppose that , for some $i\geq1,$ the numbers $m_i\mp1$ are twin primes.
 Put $n_i=m_i+3.$ Then $n_i\equiv3\pmod6$ and we have
 $$a(n_i-3)=2n_i-6$$
$$a(n_i-2)=2n_i-5,$$
$$a(n_i-1)=2n_i-4,$$
$$a(n_i)=3n_i-6,$$
We see that the main increment is $n_i-2.$ By the condition, before $m_{i+1}$ we can have only a finite set if miner increments. Suppose that, they are in the points $n_i+l_j, j=1,...,h_i.$ Then, by (7.1), we have
$$a(n_i+1)=3n_i-5,$$
$$...$$\newpage
$$a(n_i+l_1-1)=3n_i+l_1-7,$$
$$a(n_i+l_1)=3n_i+3l_1-6,$$
$$...$$
$$a(n_i+l_2-1)=3n_i+2l_1+l_2-7,$$
$$a(n_i+l_2)=3n_i+3l_2-6,$$
$$...$$
$$a(n_i+l_h-1)=3n_i+2l_{h-1}+l_h-7,$$
\begin{equation}\label{7.2}
a(n_i+l_h)=3n_i+3l_h-6,
\end{equation}
$$...$$
\begin{equation}\label{7.3}
a(n_{i+1}-3)=2n_{i+1}-6
\end{equation}
$$a(n_{i+1}-2)=2n_{i+1}-5,$$
$$a(n_{i+1}-1)=2n_{i+1}-4,$$
\begin{equation}\label{7.4}
a(n_{i+1})=3n_{i+1}-6.
\end{equation}
Note that, in every step from (7.2) up to (7.3) we  add 1 simultaneously to values of the arguments and of the right hand sides.  Thus in the fundamental point $m_{i+1}=n_{i+1}-3$ we have
$$n_i+l_h+x=n_{i+1}-3$$
and
$$3n_i+3l_h-6+x=2n_{i+1}-6$$
such that
$$n_{i+1}=2n_i+2l_h-3.$$
Now we should prove that the numbers
$$n_{i+1}-4=2n_i+2l_h-7, \enskip n_{i+1}-2=2n_i+2l_h-5$$
are twin primes. We have
$$a(n_i+l_h+t)=3n_i+3l_h-6+t,$$
\begin{equation}\label{7.5}
a(n_i+l_h+t+1)=3n_i+3l_h-5+t,
\end{equation}
where $0\leq t\leq n_i+l_h-6.$
Distinguish two case.\newline
1) Let $l_h$ be even. Then, for even values of $t$ the numbers $n_i+l_h+t+1$ are even and from equalities (7.5) we have
$$\gcd(n_i+l_h+t+1,\enskip3n_i+3l_h-6+t)=1.$$\newpage
It is easy to see that $l_h+1$ is not multiple of 3. Indeed, it is sufficient to choose $t=6.$ Thus $2l_h-1$  is not multiple of 3 and, therefore, $N=2n_i+2l_h-7$ also is not multiple of 3.\newline
Furthermore, considering $t$ not multiple of 3, from equalities (7.5) we have
$$\gcd(3n_i+3l_h+3t+3,\enskip3n_i+3l_h-6+t)=1$$
and
$$\gcd(2t+9,\enskip2n_i+2l_h-7)=1,\enskip 0\leq t\leq n_i+l_h-6,\enskip t\equiv2, \enskip4\pmod6.$$
Now in order to prove that $N$ is prime it is sufficient to use $t$ of the form $t=6u+2.$ Since $ 0\leq t\leq n_i+l_h-6,$ then $9\leq2t+9=12u+13\leq2n_i+2l_h-3$ and $0\leq u\leq (n_i+l_h-8)/6=(N-9)/12>(n_i-8)/6.$ Note that, for the considered values of $n_i(\geq183)$ we have $\frac {n_i-8} {6}>\sqrt{2n_i-16}.$ Therefore, $(N-9)/12>\sqrt{N}.$ Let $p\leq \sqrt{N}.$ Consider the congruence $12u+13\equiv0\pmod p.$ Choose a solution $u\in \{0,1,...,p-1\}.$ Then $u\leq \sqrt{N}<(N-9)/12$ and we conclude that $\gcd(N,\enskip p)=1.$ Thus $N$ is prime.

On the other hand, for odd values of $t,$ taking into account that numbers $n_i+l_h+t+1$ are odd, from  equalities (7.5) we have
$$\gcd(n_i+l_h+t-1,\enskip3n_i+3l_h-6+t)=1.$$
Note that $l_h-1$ is not multiple of 3. Indeed, it is sufficient to choose $t=3.$ Thus $2l_h-5$  is not multiple of 3 and, therefore, $M=2n_i+2l_h-5$ also is not multiple of 3.\newline 
Let now $t$ is not multiple of 3. Then
$$\gcd(3n_i+3l_h+3t-3,\enskip3n_i+3l_h-6+t)=1$$
and
$$\gcd(2t+3,\enskip2n_i+2l_h-5)=1,\enskip 0\leq t\leq n_i+l_h-6,\enskip t\equiv\pm1\pmod6.$$  In order to prove that $M$ is prime it is sufficient to use $t$ of the form $t=6u+1.$ Since $ 0\leq t\leq n_i+l_h-6,$ then $3\leq2t+3=12u+5\leq2n_i+2l_h-9$ and $0\leq u\leq (n_i+l_h-7)/6=(M-9)/12>(n_i-7)/6$ and exactly as for $N$ we obtain that $M$ is prime as well and the numbers $N$ and $M$ are twin primes.\newline
2) Let $l_h$ be odd. Then, using again equalities (7.5), by the same way, we show that the numbers $N,\enskip M$ are twin primes. This completes the induction.$\blacksquare$\newline\newline

\newpage

\indent\bfseries Acknowledgment.\mdseries \enskip The author is grateful to Daniel Berend (Ben Gurion University, Israel) for useful discussions; he also is grateful to Richard Mathar (Leiden University, Netherlands) and Konstantin
Shukhmin (Dunedin, New Zealand) for an important help in the numerical calculations.\newline\newline

\;\;\;\;\;\;\;\;


\begin{thebibliography}{1}

\bibitem 1. E. S. Rowland\slshape\enskip A natural prime-generating recurrence \upshape J.Integer Seq., v.11(2008), Article 08.2.8
\bibitem 2.  V. Shevelev, \slshape An infinite set of generators of primes based on the Roland idea and
conjectures concerning twin primes,\upshape\enskip http://www.arxiv.org/abs/0910.4676 [math. NT].
\bibitem 3.  V. Shevelev, \slshape  Generalizations of the Rowland theorem, \upshape\enskip http://www.arxiv.org/abs/0911.3491 [math. NT].
\bibitem {4}.  N.\enskip J.\enskip A.\enskip Sloane,\enskip\slshape The On-Line Encyclopedia of Integer Sequences \upshape $(http: //www.research.att.com/\sim njas)$
\end{thebibliography}
\end{document}